\documentclass[authoryear,11pt]{elsarticle}

\usepackage[a4paper, total={6in, 8in}]{geometry}

\usepackage[breaklinks]{hyperref}

\usepackage{xcolor}
\usepackage{xurl}

\usepackage{tikz}

\usepackage[greek,english]{babel}
\usepackage{amssymb}

\hypersetup{colorlinks=true}

\begin{document}

\begin{frontmatter}

\title{Integrating ethical, societal and environmental issues\\ into algorithm design courses}

\author[1]{Odile Bellenguez}

\author[2]{Nadia Brauner}

\author[3]{Christine Solnon}

\author[4]{Alexis Tsoukias}

\affiliation[1]{
    organisation={IMT Atlantique, LS2N, UMR CNRS 6004},
    postcode={44300},
    city={Nantes},
    country={France},
    }

\affiliation[2]{
    organization={Univ. Grenoble Alpes, CNRS, Grenoble INP, G-SCOP},
    postcode={38000},
    city={Grenoble},
    country={France},
    }

\affiliation[3]{
    organization={CITI, Inria, INSA Lyon},
    postcode={69621},
    city={Villeurbanne},
    country={France},
    }

\affiliation[4]{
    organisation={CNRS-LAMSADE, PSL, Université Paris Dauphine},
    country={France},
}

\begin{abstract}
    This document, intended for computer science teachers, describes a case study that puts into practice a questioning of ethical, societal and environmental issues when designing or implementing a decision support system. This study is based on a very popular application, namely road navigation software that informs users of real-time traffic conditions and suggests routes between a starting point and a destination, taking these conditions into account (such as Waze). The approach proposes to intertwine technical considerations (optimal path algorithms, data needed for location, etc.) with a broader view of the ethical, environmental and societal issues raised by the tools studied. 
     
    Based on the authors' experience conducting sessions with students over several years, this document discusses the context of such a study, suggests teaching resources for implementing it, describes ways to structure discussions, and shares scenarios in different teaching contexts. 
\end{abstract}

\end{frontmatter}

\section{Introduction}

Shortest path algorithms, such as Dijkstra's or Bellman-Ford's, are essential classics in computer science curricula. 
They also have important applications, the most emblematic of which being road navigation assistance tools which compute routes that take real-time traffic conditions into account (like Waze\footnote{Waze -- \url{https://www.waze.com}}, Google Maps \footnote{Google Maps -- \url{https://www.google.com/maps}} or Apple Plan\footnote{Apple Plan -- \url{https://www.apple.com/maps/}}). These algorithms are not "morally" neutral: along with efficiency, they bring changes that may have a significant impact on the environment and society, and this raises ethical questions that are essential to address with students as they study the technical aspects.

Our concerns arise from the fact that our students will participate in the production of digital tools which have wide impacts. Therefore, as computer science teachers (and researchers), we aim to train them to practice  a questioning of the ethical, societal and environmental issues at stake when designing or implementing an algorithm to be used for decision support purposes. The objective is to enable them to analyze the responsibility of the designers and the limits of this responsibility.

Our approach combines technical and scientific considerations when designing the tool (data needed to locate users, optimal path algorithms, state graphs, complexity, etc.) together with ethical questions, possibly accompanied by environmental and societal issues -- taking a step back on the factors and effects of decision-making and technical expertise. Beyond technical and scientific considerations, the aim is to empower the ethical vigilance of students who are trained to design automatic decision-support tools and to organize and manage their developments.

This document is intended for  computer science teachers (in artificial intelligence, algorithms or operations research, for example) who wish to open up their technical courses to societal and environmental issues. Section~\ref{sec:teaching} discusses the interests of integrating ethical, societal and environmental challenges in technical courses.

We consider a case study consisting in developing a decision support tool
for road navigation that finds the fastest route between two points for a given departure time. Additional features may be proposed by participants in order to make the tool modern, user-friendly, usable, social, economically viable, etc. The core of the case study is to take traffic conditions into account when searching for the shortest paths, and hence to modify and extend classic algorithms. The person responsible for implementing this functionality must therefore have a thorough understanding of the proofs of the appropriate shortest path algorithms -- like Dijkstra's algorithm \citep{Dijkstra59} and A* \citep{HNR68} -- in order to identify what remains valid and what needs to be adapted when choosing to complicate the problem or add functionalities. In this case, to take traffic conditions into account, the distances on the arcs are transformed into arc crossing times, which depend on the time of arrival on the arc and require adjustments of the algorithm. Section~\ref{sec:Dijkstra} provides a technical analysis of those classical algorithms and studies their extension in the case of time-dependent cost functions (complexity analysis, extension of proofs, etc.).

Section~\ref{sec:indiv-collectif}  offers some thoughts on the conflict caused by this shortest path feature between the individual interests of users and collective functioning. The aim is to analyze how a simple technical feature can cause societal and political problems. The widespread use of the tool leads to an efficient service and better data quality. However, relying on newspaper articles and on \citep{cardon2019algorithmes}, we discuss how the cumulative effect of massive use creates a conflict between the interests of the application's users and those of various groups, such as local residents, parents, truck drivers, local authorities, etc. This tension goes as far as to deny the role of local authorities in coordinating individual actions, even if countermeasures can be put in place. 

Based on the authors' experience over several years, Section~\ref{sec:reflex} presents possible concepts to structure the reflection: environmental effects classification, existing taxonomies of effects and risks with a focus on data collection analysis and ethical perspective.  Details of the contents and references that teachers can use to develop this case study in a session are given. \ref{sec:scenarisations} describes possible scenarios for implementation during lesson/practical work sessions with students.

\section{Integrating ethical, societal and environmental challenges in technical courses}\label{sec:teaching}

In this section, we first discuss the importance of addressing ethical, societal and environmental challenges in computer science students' training (Section~\ref{sec:why}). We then discuss why integrating these considerations into a technical course is relevant (Section~\ref{sec:include}). Finally, we outline some warnings about such a proposition stating the objectives and risks (Section~\ref{sec:objRisks}). 

\subsection{Why teaching ethical, societal and environmental issues?}\label{sec:why}

For several years now, the public debate has taken up the endless catalogue of problematic events following the use of algorithms, and terms such as "ethical algorithms", "ethical AI" or "ethical digital technology" are used without definition. The adjective "ethical" has become almost synonymous with every conceivable notion of justice, fairness, respect and probity, again without opening up the concepts, and with the sole identified aim of avoiding undesirable consequences. However, ethics, in the philosophical sense, is a much richer field of reflection than simply taking visible consequences into account, and allows technologies to be put into perspective, which we believe is necessary for training professionals who design or use complex digital technical systems. In order to understand the vast panorama of concepts in which a reflection of this kind can really be carried out within a technical field, we need to open up this field of reflection.

Aristotle spoke of \begin{otherlanguage}{greek} ηθική θεωρία \end{otherlanguage} (èthikè théôria) which can be understood as the rules under
which human behavior can be considered appropriate, judged or perceived as {\em good} or {\em evil} \citep{Aristotle90}. It involves looking at our actions, for which we assume we have a certain degree of freedom, and therefore thinking about this freedom, this choice, in order to move towards what would be {\em good}, in the sense of superior dispositions, in accordance with a {\em good life}. Since we are at least partly free to direct our actions, we are also responsible for them. Thinking about ethics therefore  means questioning the moral responsibility of action and how to face it by means of the major moral theories that make it possible to think about decisions, but also questioning the decision-making tools that trigger or guide action.

By taking a step-by-step approach to the design and use of socio-technical systems such as road navigation aid applications, we can establish an overview of the different areas of responsibility involved, as well as their relative limits. The first aim is therefore to enable future professionals to measure the extent to which they have a share of responsibility in the construction of these tools, which are never neutral\footnote{About the notion of neutrality in science, see for example \cite{berlan2023comment}, where author assumes that "[...] The scientific approach and attitude require the ability, at certain times and in certain spaces, to refrain from taking a position and to suspend value judgements, to silence one's convictions, in order to look at the facts without bias." (We translated) 
Since Operational Research is an operating science and not a contemplative one, the objectives in developments are never \textit{value-free}. }, and how they can equip themselves to question and work on the norms they inject into the tools, in terms of consequences, but also of underlying values and representations \citep{bellenguez2023there}. Secondly, it is essential to highlight how these responsibilities remain limited and embedded in collective injunctions and social and political contexts which form a system. The professional designer cannot bear all the responsibility for designing the standards or "rules" that are translated into the algorithms, especially as this takes place prior to use, at a time that is therefore disconnected from the action itself, and in the absence of the actors who will be the direct users of the tool (e.g. assisted drivers, who potentially have their own responsibilities), or the road users (e.g. drivers impacted by traffic even without using the tools, local residents, etc.) \citep{bellenguez25From}.
Finally, given the limits of this individual responsibility of each of the stakeholders, it is also important to understand the extended dimension of a collective (political) responsibility that professionals can help to unfold through their participation in public dialogue, by training and informing the public (institutions, service providers, distributors, users) about the issues at stake, the analyses carried out, the choices made, the technical and conceptual limits, the risks and so on.

\subsection{Why include ethical questions in a technical course?}\label{sec:include}

While some computer science programs offer courses that address these concepts of responsibility and ethics, they are generally disconnected from technical courses. By addressing both advanced technical considerations regarding tool design and ethical considerations related to its impact on society or the environment, for example, in the same practical session, the questions become more concrete. Students feel more concerned and are generally more involved as they see more directly the link with their future profession as computer engineers. A similar approach for future practitioners in industrial engineering is proposed in the map of "Limits of Industrial Thinking and Operations Research" in \cite{zarba2023ropenseeindus}.

The case studies therefore aim to address both technical and ethical questions. For instance, from a technical point of view, the case study on navigation tools consolidates the understanding of the applications and proofs of classic shortest path algorithms. They also allow us to put ethical questions into practice by considering technical choices beyond the promise they hold. The studies are mainly aimed at learners who have mastered (or are developing their skills in) the technical aspects. It allows them to understand the effects and points of view to consider when developing such a tool. As experts in the technical aspects, they are best placed to know its limitations, cost, etc.

These case studies are mainly intended for future young professionals who will have to develop decision-support tools, and it can help them clarify the issues involved in their career choices. They also help them identify the narratives underlying the organization of society in order to develop their ability to adapt to future careers.

Finally, they allow them to analyze, through practical case studies, how technical choices are politically oriented and how technical objects convey/support a vision of society. By considering the stakeholders, the effects, the consequences, and the values at stake, these studies therefore contribute to the development of an ethic of the citizen.

\subsection{The objectives and risks of such a proposal}\label{sec:objRisks}

This proposal is part of class sessions that do not follow a traditional format and may involve sharing personal values. We therefore feel it is important to include a warning that specifies the main objectives, listed below.

\begin{itemize}
     \item Analyze how and why technical tools are never neutral, despite the supposed objectivity of knowledge. Design and implementation are based on choices. 
     
     \item Discuss the underlying representations when designing 
     an algorithm, their social and political components. In particular, we can question the notion of performance~/~optimization.
    
    \item Think about what ethics means and questions, and the attention it requires.
    
    \item Take an overview, never exhaustive, to understand the risks, impacts and challenges.
\end{itemize}

This requires a careful attention to avoid "replacing" societal values that are currently being questioned and are questionable with others that will undoubtedly be questioned as well, by presenting them as "better". It is mainly a matter of engaging in a thought exercise to shift our focus to different situations and thus reveal limitations or blind spots and to show how a technical development makes sense in a given context and can lose its meaning in other situations or areas of analysis. Both teachers and students are thus faced with the absence of a perfect, immutable answer, but also with the need to work a technical approach WITH its context, which can provide a resource for analysis and a healthy limit to projections. 
Of course, this resonates or conflicts with our personal convictions and representations, which must also be questioned because they influence technical designs and the produced models.

Finally, while it is important to enable students to understand the issues at stake and their potential responsibility, it is also necessary to remind them of the limits of that responsibility. This exercise therefore also involves questioning what depends on us (as scientists, computer scientists, etc.) and what does not, in order to avoid a feeling of being overwhelmed.

\section{The case study: Technical Issues Related to Time-Dependent Shortest Paths}\label{sec:Dijkstra}

In Section~\ref{sec:comput}, we first introduce notations and recall the basic principles of shortest path algorithms, with an emphasis on the subpath optimality property which is a prerequisite for using dynamic programming. In Section~\ref{sec:conditions}, we introduce time-dependent cost functions that allow us to take traffic conditions into account, and we show that the subpath optimality property is no longer verified in the general case, thus preventing us from using dynamic programming. In Sections~\ref{sec:graphedEtat} and~\ref{sec:cond} we describe two different ways for restoring this property so that classic shortest paths algorithms may be used. 

\subsection{Computation of Shortest Paths}\label{sec:comput}

We consider a directed graph $G=(V, A)$ and a function  $c:A\rightarrow\mathbb{R}$ that associates a cost $c(uv)$ to each arc $uv\in A$. For each couple of vertices $(u,v)\in V\times V$, a $uv$-path is a sequence of vertices $P=(x_0, x_1, ...x_k)$ such that $x_0=u$, $x_k=v$, and  $x_ix_{i+1}\in A$ for each $i\in [0,k-1]$.
We note	$P[x_{i}]=(x_0, x_{1}...x_i)$ the prefix of $P$ from $x_0$ to $x_i$. The length of a path $P$, noted $|P|$, is the sum of the path arc costs, {\em i.e.}, $|P|=\sum_{i=1}^k c(x_{i-1}x_i)$. Let $s\in V$ be a given starting point. For each vertex $u\in V$, we note  $\lambda_u$ the length of the shortest $su$-path.

To compute shortest paths, Dijkstra and  Bellman-Ford algorithms use a dynamic programming principle based on the following optimal substructure property.	

{\bf Subpath optimality property}: Every prefix of a shortest path is a shortest path.

{\em Proof.\ } Let $P=(x_0, x_1,...x_k)$ be a shortest path, and $x_i$ be a vertex of $P$ such that $0< i< k$. Let $Q$ be an $x_0x_i$-path.
Let $R$ be the path obtained by replacing $P[x_i]$ with $Q$ in $P$. As $P$ is a shortest path from $x_0$ to $x_k$, and $R$ is a path from $x_0$ to $x_k$, we have $|R|=|P|-|P[x_i]|+|Q|\geq |P|$. Hence, $|Q|\geq |P[x_i]|$, and therefore $P[x_i]$ is a shortest path from $x_0$ to $x_i$. \hfill $\blacksquare$

This property allows us to recursively define the length $\lambda_v$ of a shortest $sv$-path by means of Bellman's equations:
\begin{eqnarray}
\lambda_v &=& 0 \mbox{ if }v=s \label{eq:min0}\\
\lambda_v &=& \min_{uv\in A} \lambda_u + c(uv)\mbox{ otherwise}\label{eq:min}
\end{eqnarray}
Dijkstra and Bellman-Ford algorithms are derived from these equations. More precisely, they associate with each vertex $v\in V$ a variable $d_v$ such that, at any time during the run, $d_v$ is equal to the length of the best $sv$-path found since the beginning of the run. Initially, $d_s=0$ while for any other vertex $v\in V\setminus\{s\}$, $d_v=\infty$. These variables are iteratively decreased until being minimal ({\em i.e.}, $\forall v\in V, d_v=\lambda_v$) by performing {\em arc relaxations}, where relaxing an arc $uv$ consists in updating $d_v$ as follows:
\[ d_v \leftarrow \min \{d_v, d_u + c(uv)\}\]
This instruction is a straightforward consequence of equation~(\ref{eq:min}).

Bellman-Ford's algorithm iterates arc relaxations until convergence so that each arc is usually relaxed more than once. When all arc costs have positive values, the algorithm of Dijkstra relaxes each arc at most once: at each iteration, it selects the vertex $u$ that minimises  $d_u$ and relaxes each arc outgoing from $u$ (while keeping track of the vertices that have already been treated in order to avoid selecting a same vertex more than once).

\paragraph{Possible further studies}
Dijkstra's algorithm computes all shortest paths starting from $s$, towards each vertex of the graph. When there is only one destination $f$, we may improve it by introducing a heuristic $h: V\rightarrow\mathbb{R}^+$ that sets a course for the search: at each iteration, we relax arcs outgoing from the vertex $u$ that minimises $d_u + h(u)$. The resulting algorithm, called A* \citep{HNR68}, usually allows one to find the optimal solution faster, provided that the heuristic is admissible and consistent.

It is also possible to improve performance by running two searches in parallel: a "forward" search starting from $s$ and ending at $f$, and a "backward" search starting from $f$ and ending at $s$. The optimal solution is found when the two searches converge.

Route calculators generally pre-calculate some shortest paths between a few vertices in a preprocessing step, and use these pre-computed paths to obtain better heuristics \citep{WW07}.

\subsection{Integration of traffic conditions}\label{sec:conditions}

To take traffic conditions into account, we replace the constant cost function $c:A\rightarrow \mathbb{R}^+$ with a time-dependent function $c:A\times \mathbb{R}^+\rightarrow \mathbb{R}^+$ such that $c(uv,t)$ is the travel time for arc $uv \in A$ if we start from $u$ at time $t$. Our shortest path problem then becomes a fastest path problem. Given a starting vertex $s$ and a starting time $t_0$, we denote by $\delta_u(t_0)$ the arrival time of the fastest $su$-path when starting from $s$ at time $t_0$.

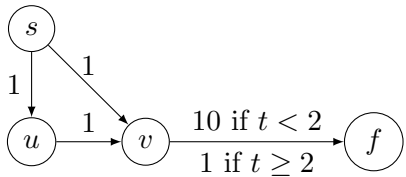
\begin{figure}[htbp]
    \centering
    \begin{tikzpicture}[scale=1.5]
	
	    \node[circle,draw] (1) at ( 0 , 1) {$s$} ;
    	\node[circle,draw] (2) at ( 0 , 0) {$u$} ;
    	
    	\node[circle,draw] (3) at ( 1 , 0) {$v$};
    	\node[circle,draw] (4) at ( 3 , 0) {$f$};
    	\draw[->,>=latex] (2) -- (3) node[above,midway]{1}; 
    	\draw[->,>=latex] (3) -- (4) node[above,midway]{$10$ if $t<2$} 
    	node[midway,below] {$1$ if $t\geq2$};    
    	\draw[->,>=latex] (1) -- (2) node[left,midway]{1}; 
    	\draw[->,>=latex] (1) -- (3) node[above,midway]{1};  
    \end{tikzpicture}
    \caption{Example for which equations (\ref{eq:min0t}) and (\ref{eq:mint}) are not valid when starting from $s$ at time $t_0=0$. Indeed, we have $\delta_v(0)=1$. Equation (\ref{eq:mint}) then gives $\delta_f(0) = \delta_v(0) + c(vf,\delta_v(0)) = 1 + 10$, whereas the path $(s,u,v,f)$ allows us to arrive at $f$ at time 3.}
    \label{fig:contre-exemple}
\end{figure}

Can Dijkstra's algorithm be adapted to solve this problem? Since all costs $c(uv,t)$ are positive, it may be tempting to answer yes. This requires a more detailed examination of the fundamentals of Dijkstra's algorithm, and in particular, encouraging students to consider whether they can extend Bellman's equations in this case. Let us consider, for example, the following equations which are a direct extension of equations (\ref{eq:min0}) and (\ref{eq:min}) to the time-dependent case:
\begin{eqnarray}
\delta_v(t_0) &=& t_0 \mbox{ if } v=s\label{eq:min0t}\\
\delta_v(t_0) &=& \min_{uv\in A} \delta_u(t_0) + c(uv,\delta_u(t_0)) \mbox{ otherwise}\label{eq:mint}
\end{eqnarray}
Are these equations valid?
The answer to this question is no, as shown in the example in Figure~\ref{fig:contre-exemple}. 
	
The problem in this example arises from the fact that the optimality property of subpaths is not verified: the fastest $sf$-path when starting at time 0 is $(s,u,v,f)$, with the arrival time at $f$ being 3, but the subpath $(s,u,v)$ is not the fastest way to go from $s$ to $v$.

Actually, the fastest path problem becomes NP-hard when the cost function depends on time \citep{ZEITZ23}. In the following two sections, we present two possibilities for restoring the subpath optimality property, thus allowing us to use Dijkstra's algorithm. The first consists of defining a new graph the size of which is generally non-polynomial compared to the initial graph but that allows us to solve the problem optimally (Section~\ref{sec:graphedEtat}). The second consists of adding a realistic condition on the cost function to guarantee the subpath optimality property on the initial graph (Section~\ref{sec:cond}).

\subsection{Exact resolution with a state-transition graph}\label{sec:graphedEtat}
	
\begin{figure}[htbp]
    \centering
	\begin{tikzpicture}[scale=2]
		
		\node[circle,draw] (1) at ( 0 , 1) {$s$,0} ;
		\node[circle,draw] (2) at ( 1 , 0.3) {$u$,1} ;
		
		\node[circle,draw] (3) at ( 1 , 1) {$v$,1};
		\node[circle,draw] (4) at ( 2 , 1) {$f$,11};
		
		\node[circle,draw] (6) at ( 2 , 0.3) {$v$,2} ;
		\node[circle,draw] (7) at ( 3 , 0.3) {$f$,3};		

		\draw[->,>=latex] (2) -- (6) node[above,midway]{1}; 
		\draw[->,>=latex] (3) -- (4) node[above,midway]{10}; 
		\draw[->,>=latex] (1) -- (2) node[left,midway]{1}; 
		\draw[->,>=latex] (1) -- (3) node[above,midway]{1}; 
		\draw[->,>=latex] (6) -- (7) node[above,midway]{1}; 
	
	\end{tikzpicture}

    \caption{State-transition graph corresponding to the graph of Figure~\ref{fig:contre-exemple} when considering only the states that may be reached from the initial state $(s,0)$.}
    \label{fig:graphedEtat}
\end{figure}
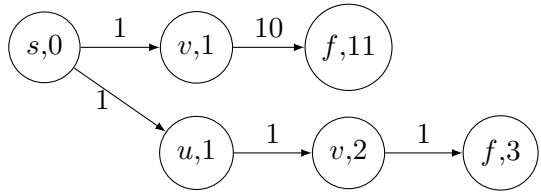

We can create a state-transition graph whose vertices (or states) are associated with pairs $(v,t)\in V\times \mathbb{R}^+$, and arcs are transitions between these states: there is an arc from $(u,t_u)$ to $(v,t_v)$ if $uv$ is an arc of $G$ and $t_v=t_u+c(uv,t_u)$, and the cost of this arc is equal to $t_v-t_u=c(uv,t_u)$. In order to reduce the number of states, we can limit ourselves to states that can be reached from the initial state $(s,t_0)$, as illustrated in Figure~\ref{fig:graphedEtat}.		

A fastest path from $s$ to $f$ when starting from $s$ at time $t_0$ corresponds to a path in the state graph going from $(s,t_0)$ to a final state $(f,t)$ such that $t$ is minimal, and this path is a shortest path in the state-transition graph. Since all costs in the state-transition graph are positive (and not time-dependent), we can use Dijkstra's algorithm to find a shortest path.
	
\begin{figure}[htbp]
    \centering
	\begin{tikzpicture}[scale=1.5]
		
		\node[circle,draw] (B) at ( 0 , 1) {$u$} ;
		\node[circle,draw] (A) at ( 0 , 0) {$s$} ;
		
		\node[circle,draw] (C) at ( 2 , 0) {$f$};
		\draw[<->,>=latex] (A) -- (B) node[left,midway]{1}; 
		\draw[->,>=latex] (A) -- (C) node[above,midway]{$2k$ if $t<k$} 
		node[midway,below] {$1$ if $t\geq k$};  
		
		
	\end{tikzpicture}
	\begin{tikzpicture}[scale=1.5]
		
		\node[circle,draw] (1) at ( 0 , 1) {$s$,0} ;
		\node[circle,draw] (2) at ( 0 , 0) {$u$,1} ;
		\node[circle,draw] (3) at ( 1 , 0) {$s$,2} ;
		\node[circle,draw] (4) at ( 1 , -1) {$u$,3} ;
		\node[circle,draw] (7) at ( 2 , -1) {$s$,4} ;
		\node[circle,draw] (5) at ( 1.5 , 1) {$f$,8};
		\node[circle,draw] (6) at ( 2.5, 0) {$f$,10};
		\node[circle,draw] (9) at ( 3.5, -1) {$f$,5};

		\draw[->,>=latex] (1) -- (2) node[left,midway]{1}; 
		\draw[->,>=latex] (2) -- (3) node[above,midway]{1}; 
		\draw[->,>=latex] (4) -- (7) node[above,midway]{1}; 
		\draw[->,>=latex] (3) -- (4) node[left,midway]{1}; 
		\draw[->,>=latex] (1) -- (5) node[above,midway]{8}; 
		\draw[->,>=latex] (3) -- (6) node[above,midway]{8}; 
		\draw[->,>=latex] (7) -- (9) node[above,midway]{1}; 
	
	\end{tikzpicture}
    \caption{Example of a graph (Left) for which the state-transition graph (Right) has ${\cal O}(k)$ states (this graph is drawn on the right for $k=4$).
    \label{fig:pseudoPoly}}
\end{figure}
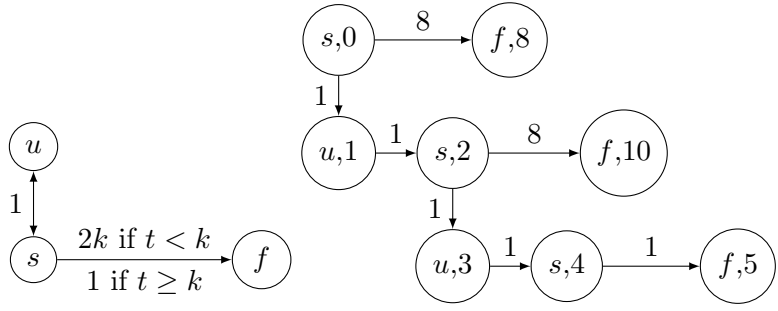

Since Dijkstra's algorithm has a polynomial-time complexity and the shortest path problem with time-dependent cost functions is $\cal NP$-hard, have we proven that ${\cal P}={\cal NP}$? No, because the state-transition graph may not have a polynomial size compared to the initial graph, as illustrated in Figure~\ref{fig:pseudoPoly}: in this example, the size of the state-transition graph is pseudo-polynomial, not polynomial. In the case where the costs have real values, it is even possible to obtain a state-transition graph of infinite size \citep{OR90}.	
	
\subsection{Addition of a condition on cost functions}\label{sec:cond}

The cost functions of the graphs in Figures~\ref{fig:contre-exemple} and~\ref{fig:pseudoPoly} are quite unusual because waiting can allow you to arrive earlier. For example, in the graph in Figure~\ref{fig:contre-exemple}, if we arrive at $v$ at time $t<2$, then it is better to wait until time $2$ before leaving for $f$, as this allows us to arrive at $f$ at time $3$ instead of $t+10$. Such a cost function is said to be non-FIFO (First In First Out) because the first to leave $v$ is not necessarily the first to arrive at $f$.

More generally, a cost function $c$ associated with an arc $uv$ respects the FIFO property if, for any pair of times $t_1,t_2\in\mathbb{R}$ such that $t_1<t_2$, then $t_1+c(uv,t_1)<t_2+c(uv,t_2)$.

Let us show that if cost functions satisfy the FIFO property, then the subpath optimality property holds. In the proof, we assume that all paths start from $x_0$ at time $t_0$. Let $P=(x_0, x_1 ... x_k)$ be a fastest $x_0x_{k}$-path. First, let us show that $P[x_{k-1}]$ is a fastest $x_0x_{k-1}$-path. Let $Q$ be an $x_0x_{k-1}$-path and $t'$ the arrival time at $x_{k-1}$ when following $Q$. Let $t$ be the arrival time of $P[x_{k-1}]$ at $x_{k-1}$.

When adding arc $x_{k-1}x_k$ at the end of path $Q$, we obtain an $x_0x_{k}$-path that arrives at $x_k$ at time $t'+c(x_{k-1}x_k,t')$. The path $P$ arrives at $x_k$ at time $t+c(x_{k-1}x_k,t)$ with $t+c(x_{k-1}x_k,t) \leq t'+c(x_{k-1}x_k,t')$ (optimality of $P$), which implies that $t\leq t'$ (FIFO property) and therefore $P[x_{k-1}]$ is a fastest $x_0x_{k-1}$-path.

By induction, any prefix of $P$ is therefore a fastest path.

Since the subpath optimality property is verified, equations (\ref{eq:min0t}) and (\ref{eq:mint}) are valid, and  we can therefore calculate a shortest path using Dijkstra's algorithm, the only modifications to be made being (i) initialising $d_s$ at the starting time $t_0$ instead of $0$, and (ii) using the following arc relaxation instruction:
\[d_v \leftarrow \min \{d_v, d_u + c(uv,d_u)\}\]
The proof of correctness of Dijkstra's algorithm remains valid in this case.

\section{Requirements analysis: When stakeholder interests are conflicting}
\label{sec:indiv-collectif}

Computer science curricula usually include software engineering courses in which students are trained to analyse requirements for designing a new application. This requirements analysis step mainly involves identifying the different stakeholders and clarifying their needs. When identifying stakeholders, students often focus on users who directly interact with the application, and the functional needs of these users are usually specified thanks to use case models (or user stories) that describe these interactions. 

For navigation applications, typical use cases are "search for a route", "compute the expected arrival time of a route", "signal an event (accident, speed cameras, etc)", or "be informed with events along the route". The main actor of these use cases is the driver or a passenger of the vehicle used for the navigation.

Besides this user, there are many other stakeholders who have various interests and who should not be forgotten during the requirements analysis step. In order to broaden the analysis step beyond the user, it is important to list as many stakeholders as possible together with examples of their interests or objectives such as, for example:
\begin{itemize}
    \item the operator, who provides the application, and wants it to be economically viable;
    \item the technical team, who designs, implements and maintains the application, and wants it to be efficient, reliable, well documented and structured, etc;
    \item economical partners (other companies than the operator), who are interested in hyperlocal advertisements, and want to have a return on investment;
    \item public authorities (at various scales from cities to countries), who design and implement public transport policies, and want the application to be consistent with these policies;
    \item local residents, who want to live in quiet, safe and unpolluted areas;
    \item the earth, whose boundaries should not be overpassed \citep{kitzmann2025planetary}.
\end{itemize}
Identifying stakeholders is important to agree on a common vocabulary, and to think about how stakeholders take into account the views of other stakeholders and are vigilant about their agency (including that of the user). In particular, participants are invited to identify conflicts of interests and analyse their causes. A goal is to understand how simple technical features can cause environmental, societal and political issues.

Numerous articles and books provide insights into the impact of digital platforms on collective functioning, \citep[e.g.][]{cardon2019algorithmes,Grumbach14, Cardon2015, Courmont21, Courmont18, bellenguez25From}.
In this section, we illustrate this on the conflicts between the individual interests of the user at the time of use and broader, longer-term  impacts on communities to which the user may or may not belong. We show that technical choices made when designing a navigation application may result in a massive transformation of traffic conditions on some roads (Section~\ref{sec:traffic}), and may negate the role of local authorities in coordinating individual actions (Section~\ref{sec:collectivite}), even if countermeasures can be put in place (Section~\ref{sec:parades}). 

\subsection{A lot of traffic in places not intended for this use}\label{sec:traffic}

In order to be effective, time-dependent shortest paths calculation must have data (traffic volume on roads, accidents, etc.), which can be collected in various ways, involving users or not. The current choice to collect data from users' devices implies that only massive use can provide a sufficient level of information. However, this massification is accompanied by an equally massive transformation of local traffic due to the stacking effect. This creates a conflict between the interests of the application's users and those of several groups, including local residents, parents, lorry drivers, local authorities, etc.

In particular, for long-distance journeys, the algorithm does not limit itself to indicating only major roads, but suggests, in order to save time (sometimes negligible), taking roads not intended for non-local traffic. The cumulative effect due to the widespread use of these applications leads to a significant increase in traffic in places not designed for this use. 

When this traffic is rerouted to residential areas, it causes nuisance and a loss of quality of life for local residents due to traffic jams, noise and pollution. This traffic, which was not anticipated when road uses were defined, can also lead to deterioration of the roads, with a decline in travel conditions and additional repair costs for local authorities and residents. It can also represent a danger when passing in front of schools or because of traffic jams hindering access to hospitals. 
These situations are widely reported in the mainstream press \citep[see e.g. \href{https://www.nytimes.com/2017/12/24/nyregion/traffic-apps-gps-neighborhoods.html}{\em Navigation apps are turning quiet neighborhoods into traffic nightmares} in The New-York Times, 2017][]{NYT}.	

\subsection{Conflict between individual interests and collective organization}\label{sec:collectivite}
	
\cite{cardon2019algorithmes} analyse the tension between individual interests and collective functioning induced by digital platforms. They show that urban regulation is shifting from a logic of collective choices that determine how cities are used to a utilitarian optimisation of platform user satisfaction. This section is partly inspired by their article.

Digital platforms are undergoing massive and widespread deployment across all areas. They are based on connecting or collaboration between users, ignoring local specificities. Algorithms therefore have the ability to regulate digital spaces, in addition to pre-existing forms of regulation such as law, the market and standards. As digital and physical spaces are intertwined, code thus has the power to influence the physical world. 

The widespread use of these tools enables more efficient service and better data quality. Thus, by developing a monopolistic tool, platforms derive their legitimacy from its efficiency. However, the user of the tool may also be the one who will later pay local taxes and who will be hindered or impeded by a particular urban development project. The criteria of the tool (immediate interest, time savings) are therefore opposed to various criteria with different horizons and groups, which are sometimes also in conflict. This opposition can raise the question of ethics and lead to a discussion on "what characterizes a good decision/a good criterion".

The interests of these user-centred platforms are opposed to those of local authorities, which represent residents living in public spaces. These local authorities derive their legitimacy from their representativeness through elections. They are responsible for coordinating individual actions through the management of public spaces by deciding on infrastructure, controlling (traffic flow, speed), maintenance, etc. Their role is not only to regulate traffic to ensure fluidity, but also to define priorities for the use of public spaces (e.g. urban transport plans).  They decide on exceptions, resource allocation and prohibitions in order to maintain a balance between populations and preserve spaces. 
The universality of algorithms hardly take these specificities into account. To integrate local authorities decisions into path planning applications, it would be necessary to work with local authorities and define a way for integrating territorial governance data, besides user traces.

The use of commercial platforms also shifts the locations or authorities responsible for decision-making and analysis, which were previously carried out by these communities. Social issues that could be discussed at the political level disappear behind technical choices made at the heart of models without consultation. Responsibilities are thus shifted, but in a more vague and diffuse way than a clear opposition between irreconcilable interests (and no user of the tool can be held individually responsible for the effects of their decisions).

The tool thus undermines the ability of public institutions to coordinate individual actions. The economic interests of platforms and their individual users negate the autonomy of these communities.
As a result, they have lost control over behaviour and can no longer organise or make decisions. The vision of public space is changing from a project for collective living to a place of rapid transit for the individual benefit of drivers \citep{agilityEffect20}. 
	
\subsection{Countermeasures and their limitations}\label{sec:parades}

Solutions to these problems exist at different levels.
Local authorities can transform road infrastructure to prohibit or discourage traffic. For example, they can  legislate to prohibit non-residents from using certain roads, restrict traffic to local residents, slow down traffic (30 km/h zones), or change the administrative categories of roads in order to cheat the algorithm. They can also make traffic more complex by using speed bumps, traffic lights, and one-way streets. Numerous press articles report on these solutions, e.g., \citep{numerama19}. These small-scale, forced urban developments impose implementation costs on residents, reduce their freedom, prevent them from taking advantage of their own knowledge, make their daily journeys more difficult, and encourage NIMBY (not in my backyard) policies, with the risk of traffic being diverted to neighbouring areas. 

Cities can also develop local transport applications, as in Lyon\footnote{Website - \url{https://mobilites.grandlyon.com/}}. However, these solutions are not adopted by passing travellers who, because of the platforms, are not encouraged to use local applications, which limits their impact and hinders their development.

At the national level, regulatory initiatives may also be proposed. For example, in France,  \href{https://www.legifrance.gouv.fr/jorf/id/JORFTEXT000046144256}{a law}\footnote{Law -- \url{https://www.legifrance.gouv.fr/jorf/id/JORFTEXT000046144256}} in 2022 requires platforms to offer more environmentally friendly routes: raising environmental awareness, displaying pollution levels, reducing journey times by at least 10\% for secondary routes, etc. This law has been reported in the mainstream press \citep{leMonde22, laDepeche22}.

The major platforms are also developing partnerships with local authorities, such as the Waze's Connected Citizens Program\footnote{\url{https://www.waze.com/discuss/t/connected-citizens-program/}}.  However, these partnerships only concern large cities. Furthermore, they create a risk of plutocracy (government by the rich): private companies in a monopoly position force dependent local authorities to create partnerships and change their policies without any guarantees of sustainability and with a loss of autonomy and digital sovereignty (e.g. data sharing). 

Taking local specificities into account can have political and social impacts. Indeed, universal and supposedly neutral  procedures can then incorporate, intentionally or unintentionally, certain intentions and ideologies. They can then exacerbate poverty and discrimination. For example, avoiding dangerous areas is discriminatory, can harm local businesses, exacerbates poverty and structural racism. When areas are dangerous for certain populations, the application can also support discriminatory policies \citep{vice2016}.

Other solutions can be discussed. For example, at the tool level, the software could provide directions to major milestones and only take actual traffic into account for the last few kilometres. From an individual perspective, users can choose to travel on foot or by public transport. From a community perspective, it is possible to limit or encourage (nudge) certain individual behaviours (carpooling, multimodal travel, shared vehicles).

\section{Other dimensions and concepts to broaden the discussion}\label{sec:reflex}

After listing the stakeholders with the students,  the discussions can be raised through different frameworks, trying to embrace the global panel of potential impacts: A classification of Environmental effects based on \citep{Horner_2016} is first analysed in Section~\ref{sec:environnement}. This opens to  classification on AI risks from \citep{AIRiskMIT} or other conceptual works which are illustrated on some issues like human agency, transparency or unsafe use in Section~\ref{sec:classification}. Then data collection is discussed in Section~\ref{sec:donnees}. These discussions lead to a reflection on how to present an ethic perspective in technical courses together with some contents for teaching in Section~\ref{sec:balancing}. Finally Section~\ref{sec:massifResponsabilite} discusses the ethics of technology in the 20th century, with a particular focus on digital technology.

\subsection{Environmental effects}\label{sec:environnement}

First of all, an example easy to introduce is the classification of the energy effects of technologies proposed by \cite{Horner_2016}. One dimension is the net energy used that can be increased or decreased by the technology and is structured as follows: ICT equipment direct consumption (embodied energy, operational energy, disposal energy), direct rebound, indirect rebound, efficiency, substitution, structural economic changes, systemic transformation. The other dimension is the scope of impact: direct, single service, complementary services, economy- and society-wide. 

The direct impact is related to the physical pressure  on the environment caused by navigation aid tools \citep{ROUSSILHE2023101296}. 
These tools depend on complex physical technical objects, which are the result of choices made by software designers.  For a 
road navigation aid tool, we can consider the following objects, for example: hardware to run the application (smartphone, battery), network infrastructure for sending data to and from data centers, servers for 
calculating routes but also storing and processing user data.
Direct effects arise, for example, during production (extraction, geopolitical issues surrounding access to rare metals, etc.), use (energy consumption during operation, etc.), and disposal (failing or insufficient recycling facilities, etc.) of this material environment.

In terms of indirect effects, road navigation tools enable better traffic management. They may therefore appear to be beneficial for the environment (better use of roads, reduction in pollution). Of course, Section~\ref{sec:indiv-collectif} points out that these tools contribute to a change in the representation and management of public space and cause inconvenience to residents and inhabitants. However, we could assume that a few people are more annoyed by route optimization that takes traffic into account, but that overall traffic is better distributed. 

This raises the question of whether the traditional and familiar Operations Research (OR) questions of system efficiency and performance are favorable to energy, environmental, and sustainability issues. At first glance, they may give the illusion that improving a system's efficiency makes it possible to address issues of sobriety and frugality through better use, an assumption that is regularly contradicted by rebound effects. The rebound effect occurs when improving the efficiency of a product/service leads to a decrease in the consumption of certain associated resources. 
The saved resources can then be used to do more of the same (direct rebound effect) or to do something else (indirect rebound effect). The expected savings in resource consumption do not occur due to behavioral adaptation \citep{gillingham2016rebound}. This effect was identified by Jevons in 1865
who observed that improving the efficiency of coal use by steam engines led to a global increase of the coal consumption.
This effect is generally divided as follows:

\begin{itemize}
    \item first order - direct rebound effect: improvements in service (cheaper, more efficient, less energy-consuming) lead to an increase in consumption of that service, which was previously constrained by saturation for example (e.g. reducing uncertainties on the risk of getting lost or traffic jams encourages people to make trips that they would otherwise have avoided);

    \item second order - indirect rebound effect: the saved resource (income, time) leads to demand for other products and services (e.g. the saved time and costs are re-invested in other activities that generate new impacts);

    \item third order - effect on the economy as a whole: macroeconomic adjustments across different sectors (e.g. the system enables autonomous vehicles and causes growth of intelligent transportation system manufacturing);

    \item fourth order - structural changes that will have their own implications: changes in human preferences and economic and social institutions (e.g. the shift towards collective behaviors -- living further, changes in tourist rentals -- leads to more traffic; Autonomous vehicles are changing the patterns in which people choose to live and work).
\end{itemize}

OR methods are directly connected to hardware choices, the devices deployed and sold, energy consumption during use, etc. However, the environmental issues around digital tools are much broader: mineral extraction and production of technical devices (terminals, servers, antennas, etc.) take place upstream of algorithmic work and are driven by economic and geopolitical actors, among others. Added to this are issues related to recycling and product life cycles, which can present specific challenges for tech players in terms of hardware (recyclability, disassembly, reconditioning, etc.) that require technical knowledge but are not algorithmic. 
Furthermore, democratization, the pursuit of innovation, and therefore massive use, amplify the impacts and escape the strictly algorithmic dimension as well. It is important not to place all issues on the same level and to understand what relates to the way systems and algorithms are designed, and therefore ultimately our ethical responsibility as computer scientists, from issues that concern digital technology, its existence or expansion, which then raise more questions about our responsibilities as citizens, innovators, and ultimately as human beings.   

\subsection{Classifications of potential impacts}\label{sec:classification}

More broadly, any digital tool may have effects to discuss. In order to enlarge discussions, we may consider some global risks assessments. For instance, \href{https://airisk.mit.edu/}{\it The AI Risk Repository}\footnote{{\it The AI Risk Repository} -- \url{https://airisk.mit.edu/}} proposed by  \cite{AIRiskMIT} analyses the effects or risks of digital technologies across two dimensions: the cause of the risk (human entity or AI, risk caused intentionally or unintentionally, and pre- or post-deployment risk)
or through one of the seven identified risk domains: Discrimination and Toxicity --  Privacy and Security -- Misinformation -- Malicious actors and Misuse -- Human-Computer Interaction -- Socioeconomic and Environmental Harms -- AI systems safety, failures and limitations.
Those domains are deeply detailed with subdomain by the authors. Each item can be used to invite students to think about some examples they may consider: which harms do they already know? how is it currently mitigated? what are the limits? which harms are not properly addressed? and so on. Here are some examples.
 
\begin{itemize}
    \item "5.2 Loss of human agency and autonomy" 
    
    This domain is described as "Humans delegating key decisions to AI systems, or AI systems making decisions that diminish human control and autonomy, potentially leading to humans feeling dis-empowered, losing the ability to shape a fulfilling life trajectory or becoming cognitively enfeebled."
    
    In this domain, one can discuss the loss of cognitive ability and deskilling:
    the use of these tools threatens the ability to read a map or navigate. \cite{aporta2005satellite} describe the risks of such a loss for the Inuit. 
    This loss of skills leads to technological dependency, which is reinforced and maintained \citep[][and \href{https://www.arte.tv/fr/videos/106608-006-A/dopamine/}{video from the series Dopanime on Waze}\footnote{Video -- \href{https://www.arte.tv/fr/videos/106608-006-A/dopamine/}{https://www.arte.tv/fr/videos/106608-006-A/dopamine/}}]{laor2022waze}.  
    This can be put into perspective  with \citep{michel:hal-04758276} which calls "to regulate the attention market and prevent algorithmic emotional governance".

    \item "7.4 Lack of transparency or interpretability" 

    Interesting questions for computer scientists can be discussed in this context: what level of interactivity is needed to explain the choices; from an algorithmic point of view, how to explain why one route was chosen over another, how to verify that the information is sincere. These issues of information reliability lead directly to the question of the openness and availability of the code (open source vs closed source).

    The transparency issue can be discussed in terms of the risks of commercial bias. Navigation applications are distributed by commercial companies which, in the case of Waze owned by Google, offer customised, targeted and localised advertising. It is conceivable that the tool could favour its advertisers (for example, by preferentially suggesting a route that passes in front of that advertiser). 

    \item "5.1 Overreliance and unsafe use"

    Even if reporting unexpected events in real time can reduce accidents, loss of skill and attention increases the danger to users (nose on the screen, commanding guidance): e.g., \href{https://france3-regions.francetvinfo.fr/provence-alpes-cote-d-azur/bouches-du-rhone/marseille/marseille-voiture-est-tombee-vieux-port-1793809.html}{Marseille: A car falls (again) into the Vieux-Port}\footnote{Article -- \url{https://france3-regions.francetvinfo.fr/provence-alpes-cote-d-azur/bouches-du-rhone/marseille/marseille-voiture-est-tombee-vieux-port-1793809.html}}. During the California fires in December 2017, users reported that Waze was directing them to dangerous areas (no reliable source, but realistic behaviour for an algorithm that directs users to areas with little traffic). For the same reasons, students report that navigation tools have already directed them onto secondary roads that had not been snow-cleared. 

\end{itemize}

On another level, the design of technical systems has been the subject of conceptual work. For example, to help students enlarging their thoughts, we can consider the concept of convivial tool from \cite{illich1973tools}, which identifies the following three desirable properties:
\begin{itemize}
    \item the tool must extend the personal range of action;
    \item the tool must not degrade personal agency by making itself indispensable;
    \item The tool must not create power relations.
\end{itemize}

\subsection{Data collection: also a technical choice}\label{sec:donnees}

Data collection raises issues related to privacy (compliance with the GDPR), possible uses (sale to brokers, manipulation of users), quality (bias), availability, etc. These data-related issues can be explored in the MIT AI Risk Repository, for example in domain "2. Privacy and security"  but also in other domains such as "4.2 Fraud, scams, and targeted manipulation".

Waze's privacy policy states that "Periodically, Waze will collect all of the phone numbers which are stored on your device’s phone contacts book." (feature "find friends"). Why are contacts needed to provide directions? This discussion raises questions about the economic model behind this software and the legal aspects of data collection, and can serve as a starting point for a discussion about the  
collected data. The aim of this section is to show that user data collection is not imposed by the functionality of the tool, but is a choice made by its designer. 

The Waze privacy policy also states that the software collects "detailed location, travel and route information". Intuitively, one assumes that software that suggests the best routes to users needs this information to locate them or determine traffic conditions. However, technical alternatives that do not collect private data can be explored. For example, here is how a GPS works:
	
\begin{itemize}
		
	\item Satellites transmit signals indicating their location and the time the signal was sent.
		
	\item The positioning device compares the signal arrival time with the transmission time. By knowing the speed at which the signal travels, it can determine its distance from the visible satellites.
		
	\item The position of the satellite and its distance determine a sphere. The device is at the intersection of the different spheres of the visible satellites.
		
\end{itemize}
	
In this system, the user does not send any information! Everything can be calculated locally on their machine.  There are also open databases that map Wi-Fi networks (e.g. {\it Mozilla Location Service} or \href{https://wigle.net/}{Wigle}\footnote{Wigle -- \url{https://wigle.net/}}). It is therefore also possible to locate oneself by scanning  the Wi-Fi networks visible to the user, again without sharing their location. 

Similarly, real-time traffic information is derived from the instantaneous speed on each section, calculated from user speeds, reported events and historical data.
Other technologies are possible. In particular, cities generally provide real-time traffic information using sensors 
such as inductive-loop traffic detectors installed in roads, for example. This information does not come from personal data and is in line with public transport policies, in the sense that it is provided for the routes that cities wish to prioritise for travel (see, for example,  \href{https://mobilites.grandlyon.com/carte}{the Lyon mobility map}\footnote{Website -- \url{https://mobilites.grandlyon.com/carte}}, which displays traffic conditions on a map of the city of Lyon).

\subsection{Balancing risks with each other is an ethical question}\label{sec:balancing}

Up to this point, we identified many potential issues such as environmental and societal impacts of the navigation tool that could lead to massive harms, which we may therefore want to address in order to reduce their impact. Nevertheless, many of those questions are difficult to address because preserving safety or avoiding harmfull effect may decrease efficiency or even lead to useless tools. In addition, people may have different opinion on what is to be done or not. In a word, those questions are relating to how to determine the right limits and many of these questionable impacts are eventually captured through a tension between human values: safety versus privacy, or justice versus freedom for example. These typically are ethical questions, in the sense that Ethics originally deals with the way to examine behaviours and actions. This term is now used also to question artefacts' actions or recommendations regarding to similar considerations. That is why such discussions may lead to consider philosophical concepts to equip reflection.

The history of philosophy is rich and cannot be summarized simply, but it is essential to emphasize that ethics is a field of reflection that concerns all (human) actions and behaviours, to evaluate how and why they are acceptable or even desirable in light of  highest human considerations. Behind that, let us consider that we are \textit{a priori} responsible, at least for a part of what we do, meaning we have free-will and we are accountable. In such a context, ethics is dedicated to question foundations of moral considerations, norms to deploy or not, and ways to make decision. This could obviously also lead to question how people who develop and deploy digital tools are acting, but a first step may focus on how ethical discussion could be transposed to have a critical analysis of the way tools themselves behave. 

All previous discussions easily lead to understand that it is never so easy to determine a good path, let alone the right path, to respect every stakeholder, every interest, to balance opposite risks, prevent any harms, have a fair process, avoid biases, and so on. A particularly difficult point for technical experts is to accept that moral philosophy will not give a definitive answer or procedure to have an "ethical tool". Despite that, it is very useful to consider conceptual frameworks to structure questioning and enlighten practices. A first step may be dedicated to discuss how any technical approach to deal with real decision or action is deeply influenced by the way the situation to address is perceived and modelled, meaning deeply linked with representations, so socially and culturally influenced. There is no "neutral" tool since a tool is designed with regard to a given goal, determined by a group of humans at a given time. Moral philosophy may help to question implicit norms, values and beliefs that are actively orienting design. 

Among the normative approaches, systems that seek to establish general reasoning to guide and justify action, we may explore: virtue ethics, deontologism, and consequentialism (of which utilitarianism is the best known). For french speaking students, we recommend, for example, the video \href{https://www.youtube.com/watch?v=BTUZ_DcBi4M}{\it
Qu'est-ce que l'éthique ?}\footnote{Qu'est-ce que l'éthique ? -- What is ethics (video) -- \url{https://www.youtube.com/watch?v=BTUZ_DcBi4M}} by Lab'Éthique or the books  \textit{L’éthique expliquée à tout le monde} ({\it Ethics Explained to Everyone}) by  \cite{droit2014ethique}, or \textit{Introduction à l’éthique  (Introduction to Ethics)} by \cite{billier2014introduction} for an introduction to these subjects.
A focus on Utilitarism could be interesting at least for two reasons: measure of utility as a way to produce a "rational" decision is often considered since the 
last two centuries, even if Herbert Simon demonstrates with his bounded rationality theory that human does not act as a \textit{homo oeconomicus}, and so limits are mandatory to explore. Secondly, optimisation methods belong to this conviction of a rational decision with optimising a given utility function, and thus can be put in perspective. J.~Bentham, father of Utilitarianism, proposed to consider as a "good" decision what maximises the overall balance between benefits and harms. Facing any digital systems with a potentially very large impact, this balance cannot just consider immediate individual interest or collective functioning (mentioned in Section \ref{sec:indiv-collectif}). One also has to embrace long-term impacts on any stakeholders (including the user that may have contradictory needs) and circles, meaning to consider an aggregation of all utilities.

Applied to  navigation applications, the utilitarian approach can be discussed as follows: The overall benefits are highlighted in many situations: reduction of traffic congestion \citep{doi:10.1177/0361198118790619}, assistance in identifying natural disasters \citep{lowrie2022evaluating}, for example, as if these advantages had to be weighed against the aforementioned disadvantages and the overall usefulness was positive. However, the reduction in traffic congestion, and therefore overall traffic, could be debated on two fronts. Firstly, the \cite{braess1968paradoxon} paradox shows that decentralised traffic can be less efficient and lead to sub-optimal solutions. Furthermore, traffic optimisation can lead to a rebound effect: reducing the uncertainties and unpleasant aspects of transport can encourage people to travel more and thus increase traffic. With regard to natural disasters, there is also a downside: applications have sometimes endangered drivers when the algorithm diverts them to roads with very little traffic because safety cannot be guaranteed (use of roads that have not been cleared of snow, roads that are dangerous due to disasters, etc.). These cases have been reported in the press, for example during the fires in California \citep{clubic2017}.

We must understand that all these approaches have limitations in practice, as they can contain conflicts, dilemmas, and blind spots when applied to real-life cases or contextual thought experiments (for our french speaking students, we recommend, for example, the videos \href{https://youtu.be/AZBDMN5wZ-8?si=sFhaf16N-WU4co6h}{\it How far will you go as a utilitarian?}\footnote{How far will you go as a utilitarian? (video) -- \url{https://youtu.be/AZBDMN5wZ-8?si=sFhaf16N-WU4co6h}} (Monsieur Phi) and  \href{https://www.youtube.com/watch?v=SbWoeR5krXs}{\it When a baby's cries become a philosophical dilemma}\footnote{When a baby's cries become a philosophical dilemma (video) -- \url{https://www.youtube.com/watch?v=SbWoeR5krXs}} (Le Monde) or the book {\it Human Kindness and the Smell of Warm Croissants} by \cite{ogien2011influenceEn}).
Normative ethics are therefore accompanied by practical ethics or applied ethics (a term that appeared in this form in the 1960s, but whose roots go back much further. Indeed, some debates on the "validity" and universalism of norms or rules could be allocated to Hippocrates, Pascal, Constant, etc.) that seek to articulate the norms or rules of normative conceptual frameworks with the realities experienced when determining action facing real circumstances, in order to derive, for example, more concrete ethical "principles"  and tools for reflection or questioning that can be mobilized within professional issues in the broadest sense. Examples include medical ethics (which dates back in part to Hippocrates but was fundamentally structured in the 20th century, \citep{beauchamp1994principles}), business ethics, bioethics, and research ethics. Such Principlism is now more and more explored in the AI field and may obviously be declined in Operational Research. 

\subsection{Massive impacts and scientists' responsibility}\label{sec:massifResponsabilite}
Then, it may be relevant to introduce into these reflections the fundamental rupture of the 20th century: ethics is brutally challenged by the tragedies designed by humans (exterminations, torture, harmful medical experiments, etc.) and by the threat of technosciences (nuclear and genetic in particular), which can give everyone the power to destroy life, transform humanity (or today to spread fake news, harass, endanger, change political regimes, pollute excessively, etc.). This gives rise to numerous reflections and feeds into applied ethics (bioethics, for example) and normative ethics on the one hand, and on the other hand leads some thinkers (in particular \cite{Jonas1979}, summarized in the video \href{https://www.youtube.com/watch?v=NPIJHZO7hQo}{The Imperative of Responsibility}\footnote {The Imperative of Responsibility (in French) -- \url{https://www.youtube.com/watch?v=NPIJHZO7hQo}}) to establish that the classical standing of ethics does not fit anymore, because the scale of human actions and power to act exceeds our ability to embrace what is at stake. In this context, reflections on the ethics of technology intersect with the philosophy of science as a whole.

Several authors have questioned the place of technology in the twentieth century (e.g. Anders, Ellul) and the changes (threats) it has caused. Stiegler's work in particular echoes Socrates' comments on writing in the Phaedrus, highlighting the inherent ambivalence of technology since the invention of writing, which was described as a \begin{otherlanguage}{greek} φάρμακον \end{otherlanguage} (pharmakon: both remedy and poison). Faced with the "digital pharmakon", capable of bringing new capabilities and unprecedented threats at the same time, Stiegler goes further than Socrates and emphasizes the need for a collective (political) framework to guide and supervise the civilization choices currently being made. This work allows students to think on the two related facets of all technology (e.g., social networks that can be used to unite virtuous movements but also to carry out massacres) that engineers must consider in their design, documentation, and precautions for use (consider the discussion proposed by \cite{coutellec2019penser}), but also on the societal dimension that goes far beyond them and requires a certain consideration for public debate, education, etc.

Close to the digital systems we are considering, the first attempts to design algorithmic systems aiming to support human decision making date back to the idea of "cybernetics" \citep{Wiener48} and the "general problem solvers" \citep{Newelletal1958,Newell1963,NewellSimon72}. At that time the critical standing with respect to such attempts has been represented by the research at the Tavistock Institute at London around the notion of "socio-technical" systems  \citep{Emery59,Beer66}. Under such a perspective we can and we should not separate the design of such systems from the analysis of their broad purpose and their long term impact upon the society. This had an obvious echo to similar topics raised within the Operational Research Societies \citep{Ackoff1974,AckoffEmery1972}. 

Later on this debate has been actualised within the Artificial Intelligence Community. In the 80s Terry Winograd published together with the former minister of the Allende government in Chile, Fernado Flores, a book partially based on the failed attempt of that government to implement a system aimed at supporting the government's policy making (Stafford Beer has been involved in that attempt). In that book, \cite{WinogradFlores86}  claim that any such system should be designed starting from the anticipation of what could go wrong instead of what is expected to work appropriately and this includes what could go wrong in the society seen as a system. 


Finally, we may address the many reflections on digital or algorithms ethics, which can be seen as a field of applied ethics. These reflections are in fact part of everything mentioned above \citep{doat2024quelle}, and much remains to be done in this currently burgeoning field \cite[see, for example, the work of][in French]{AbiteboulDowek2017, Devillers2020, Ganascia2017, Gibert2020, Grinbaum2023, HoangFourquet2024}. In addition to the legal work that needs to be introduced to students because it will change professional codes (e.g., \href{https://artificialintelligenceact.eu/}{AI Act}\footnote{European laws on AI -- \url{https://artificialintelligenceact.eu/}}), thoughts are emerging around the world that can be described as a form of extended "principlism" \citep[based on the work of][in medicine]{beauchamp1994principles}, encouraging the integration of principles of action into the training of digital actors \cite[and classified by][]{Floridi2019Unified}. The idea is to present a range of current analyses and reflections, the limitations and shortcomings of which can once again be questioned with the help of critical works \citep{munn2023uselessness} or through case studies and role-playing exercises, cultural comparisons, impact studies, etc.

Resources for teachers are available on the website \href{https://www.enseignerlethique.be/}{Teaching Ethics}\footnote{Teaching Ethics (website) -- \url{https://www.enseignerlethique.be/}} of the Catholic University of Louvain. For french speaking students, we recommend the YouTube channels \href{https://www.youtube.com/@Philoxime}{Philoxime}\footnote{Philoxime (YouTube channel) -- \url{https://www.youtube.com/@Philoxime}} and \href{https://www.youtube.com/@MonsieurPhi}{Monsieur Phi}\footnote{Monsieur Phi (YouTube channel) -- \url{https://www.youtube.com/@MonsieurPhi}} which many students are already familiar with and which offer rich cultural content on ethics.
In general, discussions on all involved aspects highlight the difficulties in balancing all the risks, but also all the potential benefits that we may wish to preserve or develop, showing how reflection on technological systems and the "good" way to design or use them is a never-ending process, which has been clearly identified as a core issue within the Operational Research community \citep{le2009ethics}.

\section{Conclusion}   

The aim of this document is to provide computer science teachers (or more broadly, science teachers) with tools to develop activities for their students that combine the technical aspects and the ethical, societal, and environmental issues of the algorithms we traditionally teach. It is based on shortest path algorithms and offers a sharing of experience and contents.

This case study is a response, tested over several years, of combining technical knowledge with the ethical questions that accompany it. As academics, it allows us to respond to external demands to integrate this dimension into our training programs, but it also reflects a repositioning in light of a broader and reshaped awareness of the moral, societal, and environmental issues surrounding techno-sciences. In particular, this study addresses the following three aspects of training engineers and academics in science and technology.

Our higher education programs, particularly in computer science, are being disrupted by the rapid evolution of socio-technical systems and their impact on the functioning of society, environmental preservation, and the questioning of human characteristics.
In order to train engineers for tomorrow's world, one solution is to offer, in addition to courses that are purely technical or purely in the  social studies, a multidisciplinary approach that combines social studies and technology, thus mixing the functioning of tools with the associated questions.

Furthermore, this study is part of a comprehensive decision-making support methodology that emphasizes the importance of teaching formal decision-making tools not as a collection of methods but as a methodology that integrates context, interest, stakeholders, etc. \citep{TSOUKIAS2008138}.
For example, the future professionals we are training will need to consider and allow for possible recourse in relation to the tools they develop or may find themselves confronted with regulations, and therefore with their legal responsibility. They will need to develop skills, for example, on how to challenge/defend a decision made by the tool in court, whether the proposed algorithm was relevant, whether its use was appropriate (see, for example, \href{https://www.cours-appel.justice.fr/paris/proces-de-lintelligence-artificielle-le-carambolage-du-siecle}{AI mock trials}\footnote{AI mock trials -- \href{https://www.cours-appel.justice.fr/paris/proces-de-lintelligence-artificielle-le-carambolage-du-siecle}{https://www.cours-appel.justice.fr/paris/proces-de-lintelligence-artificielle-le-carambolage-du-siecle}}).

Finally, this study allows us to imagine teaching beyond the transfer of technical skills, as a transfer of awareness. Our role as teachers also requires to explain the malicious acts and manipulations that can be injected into algorithms and, more generally, the power relations they create.

The proposed reflection can be developed using other examples, such as the stable marriage algorithm by Gale-Shapley (1962), which is implemented in university recruitment tools and allows to think on the concepts of meritocracy, social justice, manipulation, etc. Another example can be based on Google's Page Rank algorithm, which allows to discuss on the criteria used by search engines, advertising-based economic models, dark patterns for manipulating users, etc. We propose listing \href{https://moodle.caseine.org/course/view.php?name=ROES}{on a shared page}\footnote{Examples of educational content - - \url{https://moodle.caseine.org/course/view.php?name=ROES}} examples of educational content of this type proposed by academics who have tested them.

\ 

\textit{Acknowledgements:} The authors would like to thank the following colleagues for  passionate discussions, exchanges, references, and critical reviews: Moritz M\"uhlenthaler, Ali Yaddaden, Chloé Bonifas, Sarah Ghaffari, Marc-Antoine Pencolé and also the students who participated in this case study and enriched it.

\bibliographystyle{elsarticle-harv}
\bibliography{ref.bib}

\appendix

\section{Different scenarios}\label{sec:scenarisations}

This case study has been used in various contexts. First, we propose some content (features, impacts) developed with students to enable them to consider the possibilities and stimulate reflection. 
We describe its use in traditional courses on algorithms, operational research, programming, and computer science at universities and engineering schools, as well as in a thematic school on the ethics of algorithms.

Creating lists of features and impacts  helps to agree on the scope of the study and stimulate discussion. 
This is usually done through individual reflection followed by group work or {\em world café} style sessions.

Table~\ref{tab:fonctionnalites} (non-exhaustive list) presents some of the {features} expected of a road navigation assistance tool, possibly based on existing software. We list some proposals made during the various sessions. These features concern either the objective of finding a route, the travel environment, or social aspects related to travel. Listing the features helps to remind us, for example, how certain tools work, such as GPS location, which can be done from a machine not connected to the network, with all calculations performed locally using data sent by satellites (Section~\ref{sec:donnees}). 

\begin{table}[htbp]
    \fbox{
    \begin{minipage}{\textwidth}
    {\bf The route}
    
\begin{itemize}
    \setlength\itemsep{0pt}
    \item Locate the user
    \item Determine the shortest route between two locations
    \item Estimate arrival time
    \item Suggest alternative routes
    \item Assess the difficulty of a route
    \item Report difficulties such as traffic jams
    \item Propose the least harmful routes for external people
    \item Voice guidance
    \item Display the route on a map on screen (with different views)
    \item Estimate the cost, the amount of fuel burnt, the carbon footprint
    \item Offer alternative modes of transport (train, bicycle)
    \item Facilitate carpooling (sic.)
\end{itemize}

\textbf{The environment of the route}

\begin{itemize}
    \setlength\itemsep{0pt}
    \item Report historic monuments and points of interest
    \item Display geolocated, personalised commercial recommendations (advertising)
    \item Offer music, gamification (behaviour badges)
\end{itemize}
        
\textbf{Social aspects}
        
\begin{itemize}
    \setlength\itemsep{0pt}
    \item Share the progress of the journey with someone
    \item Enable users to call for help and be geolocated by emergency services
    \item Report difficulties (accidents, roadworks, slowdowns)
    \item Report speed cameras and law enforcement officers (always makes people smile in a course related to ethics)
\end{itemize}
\end{minipage}
}
    \caption{Features expected from a road navigation aid tool}
    \label{tab:fonctionnalites}
\end{table}

Listing the features and stakeholders helps participants identify environmental and societal impacts or ethical issues. Table~\ref{tab:enjeux} provides a few very diverse examples. Some of these can be explored in greater depth later on. 

\begin{table}[htbp]
\fbox{
    \begin{minipage}{\textwidth}
    \begin{itemize}
    \setlength\itemsep{0pt}

    \item inconvenience for local residents caused by diverting traffic onto secondary routes (due to the fact that stakeholders have been listed),
    
    \item endangering users, unsuitable route (everyone has a story to tell),
    
    \item Displaying the estimated time of arrival (may encourage drivers to drive faster), 
    
    \item reporting of police and speed cameras (one of the app's first features, which is worth discussing as it can have positive effects such as slowing drivers down),
    
    \item transmission of data such as location and speed to private companies, enabling them to infer users' habits,
    
    
    \item risk that the route may be chosen to pass by certain clients of the company (e.g. restaurants that pay),
    
    \item dependence on tools, digital prostheses/orthoses, loss of cognitive ability, reduced resilience,
    
    \item ability for users to influence the algorithm,
    
    \item monopoly situations that enable companies to impose rules on local authorities and residents, making it difficult to use local applications and therefore it becomes useless to develop them.
    
    
    
    
\end{itemize}
\end{minipage}}
    \caption{Selection of some environmental, societal and ethical issues}
    \label{tab:enjeux}
\end{table}




    
    

We now describe various scenarios from this case study that have been tested since 2020.

\subsection{Algorithmics course at INSA Lyon}

At INSA ({\em Institut National des Sciences Appliquées}, National Institute of Applied Sciences) in Lyon, the session is part of an advanced algorithmic course that takes place in the second semester for third-year students in the computer science department.

The session is divided into three parts. First, students study the problem of finding the shortest path when the costs of the arcs depend on time (Section~\ref{sec:Dijkstra}). After noting that the subpath optimality condition no longer holds, they are asked to define a state graph that can be used to solve the problem. They are then asked about the size of this new graph, the objective being to find an example where the state graph is not polynomial in size compared to the initial graph. Finally, students are led to exploit the FIFO condition on arc costs in order to restore the subpath optimality condition and thus solve the problem efficiently. 

The second step involves defining the costs of the arcs by asking a series of questions: What data can be used for this? Where can this data come from? How can its quality be assessed? How can the data be converted into a predictive model? Can this data be used for purposes other than predicting journey times?
These questions about data are not intended to address technical aspects, but rather to raise awareness of the central role of data, whether for calculating arrival times or for better targeting advertising (Section~\ref{sec:donnees}).

In the final part, students are encouraged to consider the impacts on society and the environment. After identifying the various stakeholders (and their interests) and the different technical alternatives, they are asked to consider the possible consequences for each stakeholder depending on the technical choices made (Section~\ref{sec:environnement}). 

This session fits naturally into the algorithmics course, and the technical question hook works very well with these students.


 
 


\subsection{Operational Research course in Grenoble (ENSIMAG, L3 Info)}

In Grenoble, the three-hour session takes place in an optional Operational Research course, either in the third year (Bac+3) of the Computer Science degree at Université Grenoble Alpes, or in the second year (Bac+4) at ENSIMAG ({\em École Nationale Supérieure d'Informatique et de Mathématiques Appliquées de Grenoble}, engineering school in Computer Science and Applied Mathematics). Students have already had courses on shortest path algorithms and are familiar with Dijkstra's algorithm without necessarily mastering its proof.

Students know that the session is about ethics, but they do not know the topic in advance. We use the following hook: "If your smartphone were replaced by a push-button phone (with text messaging and phone call capabilities), what would you miss the most?" Locating oneself, finding directions, or finding paths consistently ranks among the top three features considered most important.



The context is then announced: "Reflection on a navigation aid tool for car drivers". The tool must enable the following question to be answered: "A motorised user wants to find a faster route between their location and their destination". It must take into account the contexts of use, user requirements and the implementation of the tool. The students are in the {role} of the analyst who receives the order for the tool and accepts it. From this perspective, they must advise the client (who is proposing the tool) by using their methodological and algorithmic skills. 

In this context, a brief introduction to ethics is often welcome. Students then list the stakeholders (as describe in Section~\ref{sec:indiv-collectif}), the expected functionalities (Table~\ref{tab:fonctionnalites}), and the ethical issues (Table~\ref{tab:enjeux}). They are then invited to explore the issue of shortest paths in more depth, looking at the technical aspects (Section~\ref{sec:Dijkstra}) and the tension between individual interests and collective functioning (Section~\ref{sec:indiv-collectif}). Depending on the time available, the issue of data for location purposes (Section~\ref{sec:donnees}) is also discussed. The other topics presented in this article are offered as optional additions to the session. 

Students emphasise their interest in combining the technical aspects they can master, reflections on their use of technology and what they have read about it, and putting these tools into perspective. They are often frustrated at the end of the session that they have not been able to say "everything".  

\subsection{Doctoral school on ethics and algorithms}

This study was also presented at a \href{https://moodle.caseine.org/course/view.php?name=EcoleROIAEthique}{Doctoral school on ethics and algorithms}\footnote{Doctoral school on ethics and algorithms -- \url{https://moodle.caseine.org/course/view.php?name=EcoleROIAEthique}} held in Aussois (France) in May 2023.  For this thematically diverse but already initiated audience, there was no need to provide an introduction to ethics, and technical aspects were not covered in detail.  

After listing the stakeholders, features and impacts, we organised the issues into key themes. 
The participants then formed groups to work on a theme that they found interesting. During the debriefing, the participants were asked to suggest ideas for developing a travel assistance tool that did not pose the ethical problem in question. 

\end{document}